\documentclass[12pt,a4paper]{article}
\usepackage{amsfonts}
\usepackage{amsmath}
\usepackage{amsbsy}
\usepackage{amsxtra}
\usepackage{latexsym}
\usepackage{amssymb}
\usepackage{amsthm}
\usepackage{url}
\usepackage{cite}
\usepackage{pstricks,pst-node}
\usepackage{enumerate}
\makeatletter
\g@addto@macro\thesection
\makeatother

\newtheorem{theorem}{Theorem}

\newtheorem{corollary}{Corollary}
\newtheorem{proposition}{Proposition}

\newtheorem{definition}{Definition}

\DeclareMathOperator{\cone}{cone}

\DeclareMathOperator{\co}{co}

\DeclareMathOperator{\argmin}{argmin}

\DeclareMathOperator{\spa}{span}

\newcommand{\lng}{\langle}
\newcommand{\rng}{\rangle}
\newcommand{\R}{\mathbb R}

\newcommand{\mc}{\mathcal}

\newcommand{\Hyp}{\mathcal H}

\newcommand{\Conv}{\mathcal C}

\newcommand{\Kup}{\mathcal K}

\newcommand{\Aa}{\mathcal A}

\begin{document}

\sffamily

\date{}

\title{In which Banach spaces is the polar of every convex cone convex}
\author{A. B. N\'emeth\\Faculty of Mathematics and Computer Science\\Babes Bolyai University, Str. Kog\u alniceanu nr. 1-3\\RO-400084
Cluj-Napoca, Romania\\email: nemeth{\huge\_}sandor@yahoo.com}

\maketitle

\begin{abstract}

Using the notion of the normalized duality map it is characterized the
class of the uniformly convex uniformly smooth Banach spaces in which
each convex cone has a convex polar.
 
\end{abstract}

\section{Introduction}

Denote by $(X,\|.\|)$ a Banach space over the reals with its norm dual $(X^*,\|.\|^*)$. Suppose that for every nonempty closed convex set $\Conv \subset X$ the mapping  
$P_\Conv:X\to \Conv$ given by
	$$P_\Conv x= \argmin \{\|x-c \|: c\in \Conv\}$$
is well defined and single valued.	Then $P_\Conv$ is called the \emph{metric projection} of $X$ onto $\Conv$.

If $\Kup$ is a closed convex cone, then
$$ \Kup^\circ =\{x\in X:\,P_\Kup x=0\}$$
is called the \emph{polar} of the cone $\Kup$ . It is everythig a cone but in general
it is not convex .


In the case of $(X = H)$, a Hilbert space, the polar of each convex cone is also a closed convex cone. This is a fundamental concept in a vast theory with extensive applications (see, for example, the excellent monograph by E. Zarantonello \cite{Zarantonello1971}). The property of convexity of the polar is essential in this context.

As we take a step toward generalizing these ideas for Banach spaces, an important question arises, which is the focus of our note.

A primary issue lies in the opening statement of our introduction: the well-definedness of the metric projection does not generally hold in a Banach space. For instance, St. Cobzas' papers provide examples of significant Banach spaces that lack well-defined metric projections \cite{SC}, \cite{SC1}. Although a formal definition of the polar can be established in special cases, it may not correspond to the kernel of an idempotent operator defined on the entire cone.

The existence of a well-defined projection onto non-empty closed convex sets is an impportant property of
 reflexive Banach spaces. These are the spaces for which the bidual \((X^*)^*\) is linearly isomorphic to \(X\).

The uniqueness of the image of the projection characterizes strictly convex Banach spaces (see further definitions in the next section).

Thus, to ensure the desirable properties of the metric projection, the ambient space for our subsequent considerations will be 
{\bf the uniformly concex uniformly smooth reflexive Banach space $(X,\|.\|)$ with the duality pairing $\lng.\rng$.}


\begin{definition}\label{dual}
In this space is well defined and basic notion the \emph{normalized duality map}
$J:X\to X^*$ given by the relation
$$\lng Jx,x\rng=\|Jx\|_{X^*} \|x\|_X = \|x\|_X^2=\|Jx\|_{X^*}^2.$$
\end{definition}

In our ambient space $J$ is  homeomorphism with
the inverse $J^*=J^{-1}.$
(See e.g.  \cite{Drag},  \cite{Cior},  \cite{DM} , \cite{Dei}. )

This notion is the basic tool in our main theorem:

\begin{theorem}\label{ffoo}

In the uniformly convex and uniformly smooth reflexive  
Banach space $X$  each convex cone possess a convex polar if and only if 
$J^*:X^* \to  X$ maps bidimensional subspaces of $X^*$ in bidimensional subspaces of $X$.

\end{theorem}

In addition to some papers considering metric projection on cones in Banach spaces, A. Domokos and M. M. Marsh \cite{DM} focus exclusively on this question, developing the theory of variational characterization of the metric projection in Banach spaces. However, the question of polar cones is not addressed there.

Recently, A. A. Khan, D. Kong, and J. Li  \cite{KhanKongLi2025} observed that in a well-behaved three-dimensional Banach space (also considered in \cite{DM}), the polar of some convex cones is not convex.

Is the convexity of polar cones valid only in Hilbert spaces?

O. Ferreira and S. Z. Németh \cite{FN}  have exaples of  finite-dimensional Banach spaces that are not Euclidean in which every convex cone has a convex polar.

In light of the aforementioned facts, the question formulated in the title of our note becomes particularly relevant.

The problem was solved by us in \cite{NemethNemeth2025} for finite-dimensional Banach spaces, but in the form of a convex geometric problem in Euclidean space. The characterization employs the ad hoc notion introduced there of the \emph{coherent meridian} (meridians will also be used in our proofs as auxiliary tools). We shall generalize our result for Banach spaces using the standard notion of normalized duality mapping.


\section{Terminology and Preparatory Results}

To simplify the exposition, we shall assume in the sequel that $(X, \|.\|)$ is a reflexive, uniformly convex, uniformly smooth Banacu
 space over the reals, with $(X^*, \|.\|^*)$ denoting its norm dual, which is the vector space of norm-continuous linear functionals on $X$ equipped with the standard dual norm.

In our proofs, we will utilize some well-known standard notions (such as convex sets, hyperplanes, cones, etc.) and some classical results on cones and convex sets from functional analysis, as contained in the monographs \cite{Cior}, \cite{Drag}, and \cite{Dei}.


The set $\Hyp\subseteq  X$ is said \emph{a hyperplane with the normal
$a,\,a\in X^* \|a\|=1$} if 
$$\Hyp=\{x\in X:\,\lng a,x \rng=0\}.$$ A hyperplane defines two 
\emph{closed halfspaces} defined by
$$\Hyp_+=\{x\in X :\,\lng a,x \rng \geq 0\},\,\,
\Hyp_-=\{x\in X:\,\lng a,x \rng \leq 0\},$$
and two \emph{open halfspaces} defined by
$$\Hyp^+=\{x\in X:\,\lng a,x \rng > 0\},\,\,
\Hyp^-=\{x\in X:\,\lng a,x \rng < 0\}.$$

The translate of the hyperplane $\Hyp$ (halfspace $\Hyp_+,\,\,\Hyp^-$, etc.), i. e.
the sets $\Hyp +u \,\,(\Hyp_+ +u,\,\Hyp^- +u, $ etc.) are also called hyperplane (halfspace).
Sometimes we also use the notations \textit{}$\Hyp_-(a)=\Hyp_-$ (similarly for the hyperplane
and the other halfspaces) to emphasise that the normal of $\Hyp$ is $a$.

If $\Aa $ is a nonempty set and for some $x\in \Aa$ it holds $\Aa \subseteq \Hyp_- (a)+x$
then $\Hyp(a)+x$ is called the \emph{suporting hyperplane to $\Aa$ in $x\in \Aa$}.
Hence e.g. $\Hyp (a)$ is the suporting hyperplane of $\Hyp_-(a)$ in ech point of $\Hyp (a)$.

The nonempty  set  $\Conv \subseteq X$ is called \emph{convex} if
$$u,\,v\in \Conv\,\implies \,[u,v]\subseteq \Conv \textrm{, where }
[u,v]=\{tu+(1-t)v :\, t \in [0,1] \subseteq \R\}.$$

We shal also use the notation for $u,v \in X, u\not=v,\,\,[u,v =\{u+t(v-u),\,t\in \R_+\}.$
If $\Aa \subseteq X$ is a nonempty set,
$$\co \Aa$$
is he smallest  convex set containing $\Aa $ and is called
\emph{the convex hull of $\Aa$.}

The nonempty set $\Kup \subseteq X$ is called a \emph{cone} if 
 $\lambda x\in \Kup$, for all $x\in \Kup$ and $\lambda \in \R_+$.
	
$\Kup$ is called \emph{convex cone} if it satisfies

\begin{enumerate}[(i)]
	\item $\lambda x\in \Kup$, for all $x\in \Kup$ , $\lambda \in \R_+$ and if 
   	\item $x+y\in \Kup$, for all $x,y\in \Kup$. 
\end{enumerate}

A convex cone is always a cone, but  for $n\geq 2$ the two notions differ .

A convex cone $\Kup$ is called \emph{pointed} if $\Kup\cap (-\Kup)=\{0\}$. 

A convex cone is called \emph{generating} if $\Kup-\Kup=X$. 

If $\Aa \subseteq X$ is a nonempty set, then
$\cone \Aa$ is the smallest cone containing $\Aa$ and it is called 
the \emph{the conical hull of $\Aa$}.

Recal the important properties of the metric projection and the duality map in the uniformly
convex and uniformly smooth reflexive Banach space $(X,\|.\|).$

{\bf For any nonempty convex set $\Conv \subseteq X$ the metric projection $P_\Conv :X\to \Conv$ 
is well defined , single valued and continuous mapping.
The duality map $J:X\to X^*$ is vell defined homeomorphism and homogeneous mapping.}

Well known consequences of these facts are the followings:

\begin{corollary}\label{ass}
\begin{enumerate}
\item Every closed , nonemty set in $X$ is the interseclion of closed halfspaces.

\item Every closed convex cone is the intersection of closed halfspaces
 whose suporting hyperplanes contain $0$.
 
 \item For the closed convex cone $\Kup \subseteq X$ from the continuity of $P_\Kup$
 and the reprezentation $\Kup^\circ = (I-P_\Kup)(X)$ it follows that $\Kup^\circ$
 is closed  cone  (\cite{NemethNemeth2025} Theorem 1).
 
 \end{enumerate}
 \end{corollary}

\section{The polar of the hypercone}

Let $(X,\|.\|)$ uniformly convex and uniformly smooth  reflexive Banach space  $(X^*,\|.\|^*)$ its norm dual and with $\lng . \rng$
the associated duality mapping . Denote
$S$ and respectively $S^*$ the \emph{unit spheres} in these spaces, that is:
$$ S=\{x\in X:\, \|x\|=1\} \,\, \textrm{and}\,\, S^*=\{x^*\in X^*:\,\|x^*\|^*=1\}.$$

The cone $\Kup \subset X$ is called a \emph{hypercone} if $\Kup \not= X$ and there
is no cone except $X$ containing it properly. From item 2 of Corollary \ref{ass} it follows that
$$\Kup=\Hyp_-(a),\,\, a\in S^*$$
that is a hypercone is in fact a closed halfspace whose supporting hyperplane contains $0$.

\begin{proposition}\label{felt}
If $$\Hyp=\{x\in X:\,\lng a,x \rng=0\}$$    
is a hyperplane in $X$, it can be written using $J$ in the form
$$\Hyp=\{x\in X:\,\lng J(y),x \rng=0\}$$
for some $y\in S$.
Thus any hypercone  $\Hyp_-(a)$ can be written in the form $\Hyp_-(Jx)$
with some uniqe $x \in S$ and
$$ \Hyp_-(Jx)^\circ = \cone\{x\}= \{tx:\,t\geq 0\} =[0,x ,$$
that is the polar of a hypercone is a closed halfline issuing from $0$.
\end{proposition}

\begin{proof}

Since $\|a\|^*=1$, from the reflexivity and strict convexity of $X$,  by James' characterisation there exist a unique $x \in S$
with 
$$\|a\|^*=\lng a,x \rng
=\|a\|^*\|x\| $$
hence, by the definition of $J$, we must have $a=Jx$.

Take $ y\in \Hyp_-(J(x)$, that is $\lng Jx,y\rng \leq 0.$
Hence we have
$$\|x-0\|=\lng Jx,x \rng \leq  \lng Jx,x-y \rng \leq \|Jx\|^*\|x-y\|= \|x-y\|.$$

This shows that $x\in \Hyp_-(Jx)^\circ$ and by the homogeneitty of $J$ that
$$ \{tx: \,t\geq 0 \}= \cone \{x\}= [0,x  \subset \Hyp_-(Jx)^\circ.$$
From the smoothness of  $\Hyp_-(Jx)$, $S+x$ and its homothetic images with center $0$ are the only 
spheres tangent to $\Hyp_-(Jx).$
Hence $\Hyp_-(Jx)^\circ =[0,x $.
\end{proof}

\begin{proposition}\label{metsz}

Each closed convex cone  $\Kup$ is the intersection of hypercones. If

$$ \Pi _\Kup := \{x\in S: \Kup \subset \Hyp_-(Jx)\},$$
then
$$\Kup^\circ =\cup_{x\in  \Pi_\Kup} [0,x.$$
 
\end{proposition}

\begin{proof}
The first assertion is a transliteration of item 2 in Corollary \ref{ass}.

Since $\Kup \subset \Hyp_-(Jx),\,\,\forall \,\, x\in  \Pi_\Kup$,
from the fact that $^\circ$ is obviously antitonic set-mapping, it follows that
$$ \Kup^\circ \supset  \cup_{x\in  \Pi_\Kup} [0,x.$$

Suppose that $y\in \Kup^\circ$. We can suppose $\|y\|=1$. Then the hyperplane with the normal $Jy$ supports
$\Kup$ at $0$, and hence $\Kup \subset \Hyp_-(Jy)$.
This shows that $y\in \Pi_\Kup $ and completes the proof.

\end{proof}

\begin{corollary}\label{kov}

We have the relation
$$\Kup^\circ =\{x\in X:\,\lng Jx,z \rng \leq 0 \,\, \forall \,z\in \Kup\}. $$

\end{corollary}

\begin{proof} 

If $x\in \Kup^\circ ,\, \|x\|=1$,  then  by Proposition \ref{metsz} , $\Hyp_-(Jx) \supset \Kup$
and hence $\lng Jx,z \rng \leq 0,\,\, \forall \,z \in \Kup.$ 

If $x\in \{y:\,\lng Jy,z \rng \leq 0,\,\,\forall\, z\in \Kup\},\, \|x\| = 1$, then $\lng Jx,z \rng \leq 0,\,\, \forall \, z \in  \Kup$ and by Proposition
\ref{metsz}, it follows  that $x\in \Kup^\circ$.

\end{proof}


\section{The polar of the wedge}

Suppose $\mc D$  ($\mc D^*$) is a two dimensional subspace of $X$  ($X^*$). Then
$\mc M= \mc D \cap S$ ( $\mc M^*=\mc D^* \cap S^*$) will be called \emph{meridian of $S$ }
(\emph{meridian of $S^*$} ).

Consider the families of sets
$$\Delta = \{\mc M : \mc M \,\,\textrm{meridian of} \,\, S \} \,\,\textrm{and}\,\,    \Delta^* = \{\mc M^* : \mc M^* \,\,\textrm{meridian of} \,\, S^* \} .$$

If $a,\,b\in S^*, \, b\notin \{a, -a\}$  $\mc M^* (a,b) = \spa \{a,\,b\}\cap S^*$ is in $\Delta^*$. 

Denote
$$ \Delta^*_o =\{ \delta (a,b) \,\,\textrm{the meridian arc in }\,\, \mc M^*(a,b) $$
$$ a,\, b\in S^* ,  b\notin \{a,-a\} \,\,\textrm{ which does not contain diametraly oposite points}\}.$$

Take $\delta (a,b) \in \Delta^*_o.$ The convex cone
$$\mc W (a,b)= \cap \{\Hyp_-(J^*c) ;\,\, c\in \delta (a,b)\}$$
will be called \emph{the wedge engendered by $\delta (a,b) \in \Delta^*_o$}.

\begin{proposition}\label{eek}

We have the formula

$$W^\circ(a,b)=\cup\{[0,J^*c:\,\, c\in \delta(a,b)\}= \cone \{J^*(\delta (a,b))\}.$$

\end{proposition}

\begin{proof}

In wiev of Proposition \ref{metsz} we have to prove that
$$\Pi _{W(a,b)} = J^*(\delta(a,b)).$$
This relation says in fact that all the supporing hyperplanes of $W(a,b)$ are
their normals of form $J^*c $ with $c\in \delta (a,b)$. This follows from the
fact that $\delta (a,b)$ is a connected set and $J^*$ is continuous. A possible way
is to reduce the reasoning to a two dimensional plane of $X$.

Let $X_2= \spa \{J^*a,\,J^*b\}. $ Then $V:=X_2\cap W(a,b)$ is a proper cone in the plane $X_2$.
If $\Hyp$ is a supporting hyperplane to $W(a,b)$ then $\Hyp \cap X_2$ is a suporting line to
$V\subset X_2$ From the definition of $W(a,b)$ all the hyperplanes $\Hyp (J^*c),\,\, c\in \delta (a,b)$
are supporting hyperplanes for $W(a,b)$ and $d_c= \Hyp (J^*c)\cap X_2$ will be a line supporting $V\subset X_2$.
If $c$ moves continuously from $a$ to $b$ in $\delta (a,b)$ the lines $d_c$ will move
continuosly from $d_a$ to $d_b$ and hence they cover the domain$X_2\setminus V$.
Thus every supporting hyperplane to $W(a,b)$ must intersect $X_2$ in one of $d_c$
and it will be a $\Hyp (J^*c)$.

\end{proof}

\begin{corollary}\label{font}

The polar $W^\circ(a,b)$ is convex if and only if $J^*(\delta (a,b)) \subset  \mc M$ for some meridian $\mc M\in  \Delta$.
 Hence all wedges in $X$ have convex polars if and only if $J^*(\Delta^*)\subset \Delta,$ that is if $J^*$
 maps the meridians of $S^*$ in meridians of $S$.
 
 \end{corollary}


\section{Polar of the cone}

\begin{proposition}\label{veg}

The polar of every closed convex cone  $\Kup \subset X$ has convex polar if and only
this holds for every wedge $W\ \subset X$ . 

\end{proposition}

\begin{proof}
The only if part is obvious since  wedges are closed convex cones too.

For the if part take the arbitrary closed convex cone $\Kup \subset X$. We
have to show that if $x.\,y\in \Kup^\circ$, then $[x,y]\subset \Kup^\circ$.

If $x$ and $y$ are linearly dependent, then it follows that  $[0,x],\,\,[0,y]\subset \Kup^\circ$
are contained on the same line through $0$, hence $[x,y]\subset \Kup^\circ$.

Suppose that  $x$ and $y$ are linearly independent and consider $u= \frac{x}{\|x\|}  \in \Kup^\circ$ and
$v= \frac{y}{|y\|}  \in \Kup^\circ$.  Denote $a=Ju$ and $b=Jv$. From the Proposition \ref{metsz},
$\Kup \subset \Hyp_-(a) $ and $\Kup \subset \Hyp_-(b).$ Since for  $c\in \delta (a,b) $ one holds
$c = \lambda a + \mu b$ with some $\lambda ,\,\, \mu \in \R_+,$ it follows that $\Hyp (c)$ also supports $\Kup$
and $\Kup \subset \Hyp_-(c)$. Thus
$$\Kup \subset \cap \{\Hyp_-(c) : c\in \delta (a,b)\} = W(a,b).$$
Whereby
$$W^\circ(a,b) \subset \Kup^\circ.$$

Since $W^\circ$ is convex by hypothesis and $u,\,v\in W^\circ$ it follows that
$$[u,v]\subset W^\circ (a,b) \subset \Kup^\circ.$$
Since $\Kup^\circ$ is cone , we have also
$$[x,y]\subset \Kup^\circ.$$

\end{proof}


\section{Proof of the theorem}

Corollary \ref{font} and Proposition \ref{veg} give together that every closed convex cone in the 
uniformly convex and uniformly smooth Banach space $(X,\|.\|)$ has convex polar if and only if
$$J^*(\Delta^*) \subset \Delta,$$
that is if $J^*$ maps meridians of $S^*$ in meridians of  $S$,

Let be $\mc  M^*\in \Delta^*$. Then there exists some $\mc M\in  \Delta$ such that
\begin{equation}\label{egy}
J^*(\mc M^*)=\mc M.
\end{equation}

From definition
 \begin{equation}\label{ket}
 \mc M^* =\mc D^*\cap S^*,\,\, \textrm{and}\,\, \mc M= \mc D\cap S
 \end{equation}
 for the two dimensional subspaces $\mc D^*$ and $\mc D$ in $X^*$ and $X$ respectvely.
 
 Obviously
 \begin{equation}\label{ha}
 \mc D^*=\cup_{t\in \R}  t\ \mc M^*,\,\, \textrm{and}\,\, \mc D = \cup_{t\in \R} t\mc M.
 \end{equation}
 
 Since $J^*$ is homogeneous and one to one we have

 $$\cup_{t\in \R} J^*(t\mc M^*) = J^*(\cup_{t\in \R} t M^*) = \cup_{t\in \R} t\mc M.$$
 
 Hence using ( \ref{ha}) it follows that
 $$J^*(\mc D^*)= \mc D.$$
 
 Accordingly, if all the closed convex cones in $X$  have convex polars, then $J^* $ maps
 bidimensional subspaces in $X^*$ in bidimensional subspaces in $X$.
 
 Conversely, if $J^*$ maps bidimensional subspaces in $X^*$ in bidimensional subspaces in $X$,
 then since
 $$J^*(S^*)=S$$
 $J^*$ should map meridians of $S^*$ in meridians of $S$ and due to Corollary \ref{font} and
 Proposition \ref{veg} evry closed convex cone in $X$ has convex polar.
 

 \section{Cones in  the $l_p$ space}
 
 We use the fundamental results on $l_p$ spaces from the standard literature,
 see e.g. \cite{Dei}.
 
 \begin{proposition}\label{lp}
  In the Banach space $l_p ,\,p \in (1,\infty)$ every convex cone has a convex polar if and only
 if $p=2$.
 \end{proposition}
 
 \begin{proof}
 
 The sufficience is obvious since $l_2$ is a Hilbert space.
 
 To prove the necessity it suffices to consider three dimensional spaces..
 We will show that in a three dimensional subspace  $l_p^3$ 
 $J^*$ maps three linearly dependent vectors $a,\,b,\, c \in l_p^*= l_q\,\, (q=\frac{p}{p-1}$) in linearly dependent
 vectors$J^*a,\,J^*b,\,J^*c \in l_p$ if and only if $q=2$ and hence $p=2.$
 
 We have for $z\in  l_p^3$ and $a\in l_q^3$
 $$\|z\|_p=(|z_1|^p+|z_2|^p+|z_3|^p)^{\frac{1}{p}}$$
 and
 $$\|a\|_q=(|a_1|^q+ |a_2|^q+|a_3|^q)^{\frac{1}{q}},$$
and for $J^*$ the reprezentation 
$$J^*a =\frac{1}{\|a\|_q} (|a_1|^{q-2}a_1,\,|a_2|^{q-2}a_2,\,,|a_3|^{q-2}a_3).$$

Consider the linearly dependent vectors $a=(0,1,1),\,b=(1,0,1),\, c=(1,1,2) \in l_q^3$ and their image by $J^*$:
$$J^* a= \frac{1}{\|a\|_q}((0,1,1),$$
$$J^* b= \frac{1}{\|b\|_q}((1,0,1),$$
$$J^* c= \frac{1}{\|c\|_q}(1,1, 2^{q-1}).$$
Now, the above three vectors can be linearly dependent if and only if $q=2$.

 \end{proof}

 {\bf In conclusion in the $l_p$ space with $p\not= 2$ there exist wedges with nonconvex polars.}


 \section{A family of Banach spaces satisfaing the condition of Theorem 1} 
 
 (See also Theorem 5 in \cite{FN} and Proposition 5 and Corollary 1 in \cite{NemethNemeth2025}.)

Let $X$ be the $n$-dimensional real vector space ($n\geq 3$), and $\lng . \rng$ the Euclidean scalar product in $X$.

Consider the symmetric positive definite linear operator $A: X\to X$

Denote with $(Y,\|.\|_Y)$ the $n$-dimensional real Banach space  with the norm defined by
\begin{equation}\label{yy}
\|x\|_Y =(\lng Ax.x \rng)^{\frac{1}{2}}.
\end{equation}
Suppose $(Y^*,\|.\|_{Y^*})$ is the norm dual of $Y$.

Since the space is finite dimensional smooth and strictly convex, it is
 is obviously uniformly convex and uniformly smooth reflexive space.

Let  $[a,x], \,a\in Y^*,\,\, x\in Y$ the dual pairing in $(Y,\|.\|_Y)$,

Then the normalized duality mapping $J:Y\to Y^*$ is defined by

\begin{equation}\label{yyy}
[Jx,x]=\|x\|_Y^2=\|Jx\|^2_{Y^*}.
\end{equation}

Hence we have
$$\|x\|_Y=\|Jx\|_{Y^*}=(\lng Ax.x \rng)^{\frac{1}{2}}.$$

(Here in the left $x\in Y, \,Jx \in Y^*$ and in the right expression $x, Ax \in X$ and $\lng . \rng$ is the
Euclidean scalar product.)
 
 From the definition (\ref{yyy}) we must have
 $$[Jx,y]=\lng Ax,y\rng$$
 with the important mention that here $Ax\in Y^*$ and $\lng . \rng$  stand formally for the application
 of the linear continuous functional $Ax\in Y^*$ to the vector $y\in Y$
  (which formally coincides with the Euclidean scalar product of $Ax$ and $y$)
 
 Whereby we have in fact that
 $$Jx=Ax\in Y^*.$$
 
 Now $A$ acts linearly on its argument, hence so does $J$. Hence $J$ maps bidimensional subspaces 
 in $Y$ in bidimensional subspaces of $Y^*.$. It follows that $J^* $ maps bidimensional subspaces
 of $Y^*$ in bidimensional subspaces of $Y$.
 
 {\bf In conclusion: for any symmetric positive definite linear operator $A$,  in the Banach space $(Y,\|.\|_Y)$ with $\|.\|_Y$ defined by (\ref{yy}),
 every closed convex cone has convex polar.}

\end{document}